\documentclass[12pt]{article}
\usepackage{fullpage}
\usepackage{amsmath,amsfonts}
\usepackage{graphicx}
\usepackage{url}
\usepackage[hidelinks]{hyperref}

\newcommand{\GAP}{\textsf{GAP}}

\title{Building Mathematical Software -- Lessons from the Development of GAP}
\author{Alexander Hulpke\thanks{ORCID: 0000-0002-5210-6283}\\[1ex]
\small Department of Mathematics, Colorado State University\\
\small 1874 Campus Delivery, Fort Collins, CO 80523-1874, USA\\
\small \texttt{hulpke@colostate.edu}}
\date{}

\begin{document}

\maketitle

\begin{abstract}
Mathematical software has particular needs that are not all covered by the
existing literature on software engineering or algorithm design. This paper
looks at some of these aspects, in the concrete context of the computer
algebra system {\GAP}. It also studies the potential impact AI systems could
have on the development process for mathematical software.

This is the submitted version of the paper, whose official accepted
version is at
\url{https://doi.org/10.1365/s13291-026-00317-0}
\end{abstract}

Building a large software project, in particular one that survives for decades,
is not necessarily part of the standard job description, or training, of
mathematicians. It raises questions, not all of which are of mathematical
nature, or have mathematical answers that are helpful (e.g. on the question of
whether one should consider implicit embeddings of domains: Is $2$ an
element of the field $\mathbb{Q}[x]/(x^2-3)$? Is $1$ considered to be equal to the
polynomial $1\cdot x^0$?).

There are also questions about group collaboration, in particular in the
context of free software, or the technicalities of software engineering.
These aspects might be unfamiliar to mathematicians, but they are already
covered extensively by existing literature in management and computer
science \cite{brooksmanmonth,raymond1999cathedral,maguire1994debugging}.

But there are other questions and problems that are connected naturally to
the unique aspects of mathematical software. The developers of many systems will have
encountered them, but they typically have not been disseminated outside the
immediate group of system authors.

The goal of this paper is to make a small step in this direction. It
describes some aspects of software development that the
developers of {\GAP} encountered over time, and which (at least in my view)
are surprising.
The solutions and approaches presented here are
certainly collective efforts, while the way they are presented will have
been tinged by the author's personal view (and passing time). Any
intellectual credit for them should not accrue to the author
of this article, but to the authors of {\GAP} as a group, foremost Martin
Sch\"onert, to whom much of the fundamental system design is due.

As Section~\ref{vibecode} will show, many of these observations
seem to remain relevant even in times of artificial intelligence
and coding agents.

This collection of items is, by its very nature, a disparate one. There are
mutual connections amongst the topics, but these do not necessarily follow
the linear form implicit in an article.

\section{Robustness and speed as necessities}

Successful academic software will find many users in academia.
But there are more.
Mathematicians train many more students than there are positions in academic
departments. It thus will not come as a surprise that there are people
outside academia who have a solid mathematical training and are interested
in trying mathematical software that was originally built just for an
academic setting. Indeed, this is possible with mathematical software in a
way it would not be with, say, a particle accelerator.
Mathematical
software thus is one of the few ways in which academic mathematics on the
research level will interface directly with a broader population and will
find users outside the direct academic setting -- a situation that I consider
very positive.

But such users often have a different computer setup than is typical in an
academic setting; furthermore, they will typically not have a system
administrator whom they can ask to install the software. But they will use
the system in a nontrivial way. Documentation and installation instructions
thus need to consider a much broader pool of users than the program authors
might have envisioned initially.

Such a broader audience also means that software will encounter inputs that
none of its authors had foreseen, and this includes the most basic ones.

\subsection{Ensuring robustness}

Small cases are important. (The well-known false induction proof that all
horses have the same color fails precisely at the step from 1 to 2.) In
the context of mathematical software, the base cases are the boring
situations in which the task is trivial and not much has to be done: the
trivial group, the 0-dimensional vector space. But the code written for the
general case often stumbles on these basic situations. The generator list
might be empty. There is no element of prime order. The group has no orbit
of length bigger than one on the set of natural numbers. These properties
make many existing routines fail, by accessing elements of an empty list, or
searching forever for elements that do not exist. But of course these cases
will arise, not just in theory, but in practical calculations, indeed as
base cases when reducing problems to smaller ones. Initial releases of
{\GAP} had an
embarrassing series of bugs arising in small examples. Ultimately, we
resolved this by writing a ``sanity check'' routine that goes through all
groups of small order (say, up to order 100), and just calculates everything
(subgroups, character table, conjugacy classes, automorphisms) we could
think of. This is not only to check results for correctness, but even just
to make sure that the calculations run through. Initially, this check found
a large number of base case errors; fixing them of course made {\GAP} far
more robust.

A related lesson is the general necessity of tests, but the burden of
running them regularly (who does it? On what machine? How do we check
completion or report problems?) delayed their introduction for a long time.
(Checks were only run, by hand, when preparing a release.) Only the
introduction of web-based automated testing (we now use
GitHub workflows) made successful tests a requirement for any integration
of new code.

We are now running two classes of tests: one consists of tests written by
the code authors, covering
what they consider a representative set of examples. The second is simply a
collection of errors reported in the past by users. This is the closest
approximation we have of what users want to run; it also avoids the
embarrassing situation in which the same (or a related) error arises again after
it has been fixed.

One of the main tools for tests is the use of assertions.
They are checks,
built into the code, that verify conditions that the code assumes will
always be satisfied (say, the program calculating a correct result, or
different data being consistent), but which in practice might fail because
of a bug. These tests can be turned on or off globally. The code thus can
verify, on demand when testing or debugging, that all these conditions are
satisfied, but the tests do not impose a runtime cost in normal use, as they
are turned off by default.
Such checks can be indispensable in debugging
code that was written long ago, and in ensuring the correctness of code as
subroutines get replaced with newer versions. But as layers of code build
on top of
each other, the cost of such checks becomes harder to control. Consider,
for example, a routine constructing a homomorphism by providing source group
generators and their images. It is entirely plausible (and {\GAP} does this) to
verify the validity of these parameters, e.g. by computing a presentation of the
group and checking that the generator images satisfy the relations. But such a test is somewhat costly, which becomes an obstacle if
there is a higher-level routine that constructs a large set of
homomorphisms. If each of these constructions implies a validity test,
the overall cost will slow any test down so much that it becomes
infeasible. We thus ended up agreeing on a hierarchy of assertions, some of
which would always be part of standard tests (namely those that do not incur
much of a time cost), while others would only be turned on in particular
cases if a more thorough verification is needed of a particular routine.

\subsection{Making algorithms run fast}

The speed of software (and the associated culture of
``benchmarks'') has been a field of competition from the earliest days of
computing. But for mathematical software there is a deeper reason for
needing speed. Mathematical problems are often not bounded in size (or
at least not by any realistic bound: even the study of the finitely many sporadic groups
requires orders up to that of the Monster group, roughly $10^{54}$). Moreover, mathematicians often like to push
computation as far as is feasible, and such a strategy is clearly limited
by the efficiency of implementations.

This means that na{\"\i}ve basic algorithms (say, using element lists or
multiplication tables) do not really serve the user community.
Implementors need to keep efficiency in mind, e.g. by using
effective representations and considering not only asymptotic cost, but also
concrete real-world performance of the implementation on research-level
input. Some of these aspects are discussed in published descriptions of
algorithms, but others date back to ideas that were only ever discussed
informally in the research community. For example, while variables in {\GAP} are
typically pointers to objects, it is desirable, for performance
reasons, to allow for quick operations on small integers.

In {\GAP} (and presumably also some other systems) this is implemented in
the following way: A {\GAP} variable typically is a pointer to (i.e. the memory
address of) the object. Such a pointer is word-aligned (so the address is at least divisible by 4).
A pointer thus always ends in bits $00$. We thus are able to store other small
objects in the same memory space, distinguishing them from pointers by
setting the two lowest bits differently. The most prominent application of
this trick is for
small integers. A ``small'' integer $a$ with $|a|< 2^{60}$ (on a
64-bit system) is stored as
$4a+1$ in memory and thus has lowest bits $01$.  If we perform -- say --
addition on two such objects $4a+1$ and $4b+1$, we calculate their
sum $4(a+b)+2$ {\em directly on the processor, without going into a general
object addition routine}, and reset the two lowest bits back to $01$. This
makes for much faster addition, and we only lose 3 bits (2 for the $01$
indicator and one for detecting integer overflow).

The performance impact can be seen easily in the following experiment,
in which we
perform a large number of integer calculations (a matrix product)
and scale up the numbers. Times are in milliseconds:
\begin{verbatim}
gap> m:=List([1..200],x->List([1..200],y->Random(Integers)));; # seed matrix
gap> powt:=function(a)local t;GASMAN("collect");
>     t:=Runtime();a:=a*a;return Runtime()-t;end;;
gap> List([0..50],x->powt(m*2^x));
[ 5, 6, 6, 6, 6, 5, 6, 6, 6, 6, 6, 6, 6, 6, 6, 6, 6, 6, 7, 6, 6, 6, 6, 6, 6,
  6, 5, 35, 101, 189, 252, 272, 201, 204, 206, 209, 208, 206, 198, 203, 199,
  205, 202, 206, 205, 206, 205, 204, 203, 210, 200 ]
\end{verbatim}
(The \texttt{GASMAN} command minimizes the risk that
garbage collection interferes with timing.)

The output shows a very notable jump in timings for $x=27,\ldots,30$
(values 35, 101, 189, 252). This is
explained by the fact that the matrix entries lie between $-10$ and $10$,
and the absolute values of the entries in $m\cdot m$ are mostly
smaller than $2^7$. Multiplication (twice) by $2^x$ thus will, for $x$
beyond $(60-7)/2=26.5$ push these entries out of the range for small
integers. This quickly extends to intermediate results, explaining the
further growth, until for factor $2^{30}$ basically the whole calculation
requires long integers.
\medskip

Robustness and speed are inherent qualities that need
to be considered while building the system, not just polish that could be
added at the end.

\section{Implicit assumptions and general usability}

As mathematicians, we pride ourselves on the exactness and specificity of
our formalisms and definitions. But when we work, often a certain sloppiness
comes in, based on the fact that we most often work on a specific set of
structures. In a paper this is easily done: there will be a sentence at the
start stating that all groups considered are finite (or solvable), or that
all rings have a one. But in a program this is not that easily done, and
is indeed ill-advised if we want to serve a broad set of users. For example,
in my own research background, rings arise mostly as sets over which group
representations are written. They often come from finite fields, and an
assumption about the ring having a one is often implicit. (Of course this
property already fails for ideals in the group algebra. But we use a
different part of ring theory for it than we do for the matrix
coefficients.)

In a software system this kind of sloppiness will run into problems. There
are perfectly legitimate reasons to study rings without one, and researchers
in these areas also would like to use existing software. In some cases such
an issue is easily fixed by adding case distinctions to existing code.
In other situations it became necessary to rethink concepts from the
start to handle cases in which the original author assumed a more
restrictive structure than that used by the broader community.

{\GAP}, in general, tries to be careful in not mixing up domains
-- the integer $1$ is different from the polynomial $x^0$, and
$2$ is not an element of the quotient ring $\mathbb{Q}[x]/(x^2-3)$, but only
of the isomorphic cyclotomic subfield $\mathbb{Q}[\sqrt{3}]$.

Once one starts with considering such distinctions, one encounters many situations that
had initially been considered completely solved, and as standard textbook
material, which in fact are much trickier. An easy example would be to
consider the definition of a GCD in the ring of even integers, in which $6$
and $4$ have no GCD.  Nor is this
something one could avoid by more careful generalization in advance.
If we wanted to do this, we would be reworking \emph{Principia Mathematica}
\cite{whitehead1910principia} and (as in that book) would barely ever arrive at $1+1=2$.

\section{Types}

Every object in {\GAP} has a type, as part of {\GAP}'s language.
This type represents {\em current} system
knowledge about the object -- what can be done with objects, what
properties (or derived information, such as the conjugacy classes for a
group) are known, how it is stored internally -- rather than mathematical truth.

The purpose of the type system is -- despite the temptation of formalization it has posed for
some users in the past
-- to allow
for the selection of the best routines for particular tasks, rather than to
allow formalization of mathematics (or maybe even formal proofs, as e.g.
the proof checker Lean \cite{demoura2021lean4} does).

While the type system uses terms (such as ``category'' or ``family'' --
see~\cite{breuerlinton98} for definitions) that
are also used in formal mathematics, there is not a clean correspondence
between these {\GAP} terms and the mathematical concepts; also logical
implications will not always work as a mathematically trained user,
inexperienced with the system, might expect:

Types represent the current knowledge, and thus might change if more
information is found out about an object. But {\GAP} will not try to find
out such information on its own, just to make the type as precise as possible.

What {\GAP} records as capabilities for algorithms is not always mathematical
truth but sometimes the analog of a medieval class system in which
abilities are fixed at creation: Any square matrix can potentially be
inverted (though the calculation might fail) while an element of a finitely
presented semigroup can never be inverted (even if an inverse exists in the
semigroup).

There are some automatic implications, and inheritance of properties (e.g.
if $M$ is an $R$-module and $R$ is a field, then $M$ is a free module).
But not every mathematical implication can be coded -- a finitely
presented semigroup will never be considered a group, even if all elements
have an inverse.

\subsection{Using types to solve problems}

The types in {\GAP} are used to select, for a particular task, the best
routine (called {\em method}) available.

For example, determining the conjugacy classes in a finite group can employ better
methods if the group is known to be solvable, known to be simple
(in which case a random search performs well),
or known to be a direct product.  (There
is in fact a small subtlety here: It is not sufficient for the group to be known to be a
direct product; we need to know the decomposition into direct factors.)

Crucially, only the given information is used; {\GAP} will not (unless an
algorithm explicitly asks for it) determine the values of unknown
properties just so that a better method could potentially become available.
\smallskip

On a technical level (see \cite{breuerlinton98} for details), the selection
of the best available method
is done through bit lists, carried with the types, that indicate the information
known by the system about the object.

When a method for an operation -- say compute the order of a permutation
group -- is made available (``installed'') to the system, it is done by
declaring
which bits must be set, for which arguments, so that the method can apply.
It is possible that there are several methods installed that could apply.
In such a situation, the method which is installed with the largest number
of requirements (which typically is the most restrictive one) is chosen,
assuming that it only exists because it is better than those with fewer
requirements. (There are extra rules for breaking ties or for indicating particularly
valuable properties.)

This model (and in particular the policy of not computing extra information
automatically, even if it could help) can sometimes lead to surprising
situations. For example, for a newly created finitely presented group, {\GAP} might
refuse to do most calculations, only to perform them happily once assured
(or explicitly asked to verify) that the group is finite.
But {\GAP} will not attempt to test this on its own.
\smallskip

The underlying reason for this behavior is that the cost of a test might be
substantial, and thus tests should only be done if explicitly requested by the
user.

As a rule of thumb, for groups, finiteness and knowledge of the group order,
followed by solvability, are the most useful properties for the selection of
better methods.
\medskip

But the cost of algorithms depends not only on mathematical properties
of the group, but also on the way it is represented on the computer. (For
example, it seems plausible that most calculations in $S_6$ are done more
naturally in degree 6 than in a larger degree.) But finding such a better
representation (if it exists) will take time, as will the evaluation of an
isomorphism between two representations.
Thus, again, {\GAP} typically will not attempt to do this,
with one major exception:

For some classes of groups, the best available methods in {\GAP} are -- due
to a lack of other implemented methods -- actually to work in an isomorphic
image, the isomorphism being induced by a suitable group action. Two
prominent cases of this are matrix groups over finite fields (for which a permutation
representation on the vectors of the underlying vector space can be taken) and
automorphism groups (which act on the underlying group). For these groups, a
flag is set indicating that all calculations will be routed through
a special homomorphism, called \texttt{NiceMonomorphism}.

It took altogether several years to put {\GAP}'s type system into place, but
overall it has aged impressively well. Over time, there were minor adjustments, but the
structure still stands basically the same 30 years later and works well in
practice. Credit for it needs to go to Martin Sch\"onert and
to the authors of
\cite{breuerlinton98}.

\subsection{Types, the computer science view}

Alongside these mathematical types, there is also the concept of data types
in computer science, most prominently
in the idea  \cite{wirth1976algorithms} that algorithms
inherently are built around data types. Here the situation in {\GAP}
is a bit different.

Fixed types make it easy to determine storage requirements, and allow for
direct access to an object's data, avoiding an intermediate access function. But
types can also be a nuisance, allowing a function to be written only for
particular objects. (When, a long time ago, I learned Pascal, I was
infuriated that I would only be able to write a sorting procedure for lists
of a particular length -- the data type of the parameter -- and would have
to basically duplicate the function if I wanted to sort lists of a different
length.)
\smallskip

By not being strictly typed, {\GAP} avoids all this rigidity and
inconvenience in its language; indeed it is possible to write very generic
functions that accept all kinds of arguments. On the user level, lists,
sets, vectors, and strings are all the same data type, namely lists.
This avoids many type conversions or special functions. But this freedom
comes initially at a high cost to performance.

As an aside, it can be
interesting to compare this typing philosophy to what is done in the system
Magma~\cite{magma} (which provides similar functionality):
Magma is much more strictly typed, so a permutation in $S_5$ cannot
automatically be multiplied with a permutation from $S_{10}$. As a result, Magma
code contains a lot of type coercions. But these fixed types make
it possible for many routines to be inherently faster, since they can make
strict assumptions about what the objects are.

What we therefore found,
over time, was that {\GAP}'s objects, internally, had to have a more
specified type that goes beyond what {\GAP}'s language determines. This
assured certain properties -- e.g. being a list of finite field elements --
so that algorithms could assume that certain more efficient methods could be
used, e.g.~binary search in a sorted list. These properties are stored
through internal flags of the system that sometimes can be checked but
cannot be set explicitly by the user. This includes whether a list is
sorted, whether the entries of a list are all of the same type, or whether a
list can be changed. The latter seems to be a weird property at first glance
-- why does this matter?

Consider a list of vectors, remembering that
vectors are lists themselves. If the entries of the vectors could be
changed, we could not store that a list of vectors is sorted, and we would therefore
lose out on efficient binary search. {\GAP} thus allows one to force objects to
be immutable, meaning that changing their entries is explicitly forbidden.
Lists that contain immutable objects may store and remember that they are
sorted.

Another example is strings. Strings, as lists of characters, will be printed
differently than usual: we want them to be printed as \texttt{"ABC"} rather
than \verb+['A', 'B', 'C']+. Thus there is an internal flag that indicates whether a list is actually
a string. But what happens with the base case? I remember a discussion we
had, probably around 1995, about whether the empty list that does not contain
integers should be exactly the same data type as an empty list that does not
contain characters. The answer to this esoteric question in fact was that
the empty list is an object of its own type and changes type as soon as
entries are assigned.

I do not claim that this soft typing with internal checks is
the optimal solution. The fixed typing in Magma certainly has advantages,
but makes writing code more tedious. I suspect that the correct answer
actually lies not in the system design but between the human and the
machine. A preprocessor (nowadays maybe one with some artificial
intelligence built in) could take the user-supplied code, without fixed type
assumptions, and insert type conversions wherever necessary. This would be
the best of both worlds. But it probably is too late to revise such a design
decision without a major rewrite of the system.
\medskip

Related to this flexibility in typing is the very
disappointing performance of the {\GAP}-to-C compiler. Na{\"\i}vely,
one would have expected such a compiler to provide a speedup
of a factor of 10 or even more. But because variables are not typed, such a
compiler (at least a compiler using technologies available in 1997) would
not be able to automatically use the most efficient routines for -- say --
addition of small integers, but every operation, even a simple \texttt{x:=x+1},
would still have to be fed through the generic mechanism
that decides at runtime which routine to call for an operation. This overhead cancelled out any
speedup that would have been obtained by writing in C. Instead, the best use
of this compiler turned out to be as a starting point for human-written C
code, providing clear guidance on how to access the innards of a {\GAP} object
and how to interface a C function to the system. The user would then
fetch the data from {\GAP} once, put it into typed C variables, and avoid operation
dispatch in favor of direct low-level calls, say for small integers. This
hand-edited code then indeed shows the desired speedup.
\smallskip

In all these aspects, types are there to select effective algorithms, rather
than to represent mathematical truth or to satisfy computer science
orthodoxy.

\section{Memory, but not too much}

In the same way that {\GAP}'s language does not require the user to force
a type on variables, an underlying garbage collector (a tool that reclaims memory that had been
occupied by objects that are no longer in use) allows objects to grow in
size and take on new parts without the user ever having to worry. Built-in memory
management and garbage collection take care of it all.

But (this is again the problem of mathematicians pushing computations
as far as they will go) the question of how much memory to allocate actually
is not trivial. On a
multi-user system, or -- nowadays -- on a computer that also holds the
user's e-mail, text editor, and web browser, allocating memory until none
is available will cause problems for other processes. This turned out to be
particularly problematic because of one feature of {\GAP}, namely the
automatic fallback
to generic methods if nothing better is available.
If you ask {\GAP} a question, it will always try to answer it. But if there
is no good method
available for answering it (or the object does not know sufficient
properties that would allow {\GAP} to use a clever method), then {\GAP} will
default to a very generic method that relies on computing explicit element
lists and which will quickly allocate huge amounts of memory.

After several incidents (the author does not deny that he has been the
source of some of them), in which a carelessly growing {\GAP} job
inconvenienced other users, it was decided to put a hard limit on {\GAP}'s
memory allocation. But that in turn created problems if the user did not
anticipate that they would require a bit more memory at a later stage in the
calculation. The result is the current, practically usable but not very
elegant compromise: if {\GAP} allocates memory (and it will always allocate
more than what is minimally needed, to avoid too many memory expansion
steps) beyond the initially set limit, {\GAP} will act at the next
opportunity
as if the user had interrupted the calculation, prompting the user to explicitly
continue the calculation, which also gives the user the option to stop a
runaway memory-grabbing job.

\section{Groups via generators}
\label{secgens}

Groups in {\GAP} are typically given by a generating set, taken from a larger
group (for example, a symmetric or general linear group), or as words
in generators which are subject to relations. Every element of the group can
be represented as a finite word in the generators, and homomorphisms can be
described through group generators and their images under the homomorphism.

This is not how groups are typically described in introductory textbooks,
though it is in fact much closer to what happens in linear algebra.
For groups defined in such a way, a principal tool, the orbit/stabilizer algorithm,
creates problems as one proceeds to subgroups.

\begin{figure}[t]
\begin{center}
\includegraphics[width=8cm]{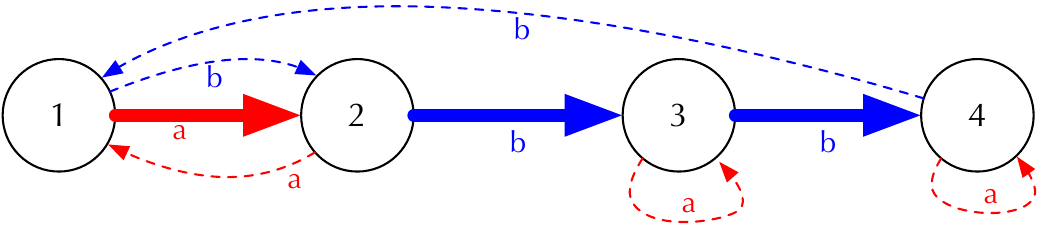}
\end{center}
\caption{Schreier graph for the orbit of 1 under $S_4$}
\label{schreiergraph}
\end{figure}

Let me briefly describe, in an example (see~\cite{holtbook}
or~\cite{hulpkebook} for a detailed description of the algorithm),
what is happening:
We consider the orbit of the point $1$ under the group $G=S_4$, generated
by $a=(1,2)$ and $b=(1,2,3,4)$, see Figure~\ref{schreiergraph}. (We use
{\GAP}'s convention that groups act on the right, $2^a$ is the image of $2$
under $a$, and permutations multiply left to right.)
We start with the point $1$ and repeatedly apply all generators, and their
inverses (in the case of a finite group we can skip inverses, as they will
also be powers with positive exponent),
to all points. The bold arrows indicate
discoveries of new points; the dashed arrows lead to points that were
discovered before.
The bold arrows form a spanning tree of the graph. Forming products along
its edges, including the identity for 1, gives us {\em representatives} for the points (and these form
a transversal for the cosets of $\mbox{Stab}_G(1)$ in $G$):
\[
1=1^{\mbox{\scriptsize id}}, 2=1^a, 3=1^{ab}, 4=1^{abb}.
\]
In general, the cost of computing such a transversal
is the total number of edges, that is,
the number of generators times the orbit length.

For an orbit of length $n$ and $m$ generators, there are $n\cdot m$ arrows
in total, the spanning tree contains $n-1$ arrows. So there will be $n\cdot
m-(n-1)=n(m-1)+1$ dashed arrows.

The dashed arrows do not give us new points, but if we combine a dashed
arrow with the spanning tree paths for its start and (the inverse of the path for its) end, we get products
that move 1 to itself: $ba^{-1}=(2,3,4)$, $a^2=()$, $ab\cdot a\cdot
(ab)^{-1}=(2,4)$,
$ab^2a(ab^2)^{-1}=(3,4)$, $ab^3=(2,4,3)$.
These elements are called {\em Schreier
generators}; Schreier's theorem shows that together they generate the
stabilizer of $1$ in $G$.

There is typically redundancy amongst the Schreier generators, but no
general result can predict it; for a free group all are needed
(Nielsen--Schreier theorem). Redundancy
has to be detected case-by-case.

The number of dashed edges in the
graph is large.  For example, if we
have two generators and an orbit of length $1000$, there will be
$1000(2-1)+1=1001$ dashed edges, leading to potentially $1001$ Schreier generators for
the subgroup. Elimination of the identity and duplicates can remove some,
but that will be far from optimal.

This creates
problems if we now, iteratively, want to compute a new orbit under the
stabilizer.

This problem was already observed by Sims \cite{sims70} in his algorithm to
compute the order of a permutation group through iterative orbit stabilizer
calculations (if the last stabilizer is trivial, the group order is the
product of orbit lengths): If $x\in G$ and $r$ is the fixed
transversal representative with $1^x=1^r$, then $xr^{-1}\in\mbox{Stab}_G(1)$.
Knowing such orbit structures not only for the group, but also for a point
stabilizer, and further stabilizers therein (a {\em stabilizer chain}) thus
allows for membership testing. One can then test each new Schreier generator
obtained for redundancy through a membership test in the stabilizer.
But the large number of tests required makes this expensive.
\smallskip

Similar problems arise in many other situations in {\GAP}, namely whenever a
subgroup is created with a larger number of generators than would be
minimally possible; often there are redundant generators involved. The
reason for this can indeed be an orbit algorithm; it could be that we are
taking the preimage under a homomorphism and thus combining generators for the
kernel with preimages of generators in the image; or it could be the result
of building up a subgroup iteratively from smaller ones.

Again, this means that
orbit algorithms for this subgroup will be unnecessarily expensive. In some
cases we could do a redundancy test, but in other situations (if the group
in which we are working does not have good methods for this or is very large)
this might not be feasible. And redundancy tests cost time on their own. At
present, the process is very heuristic: if we know that there will be a
potentially expensive orbit calculation, we force a generator reduction, but
it also seems wasteful to do such a reduction every time.

The answer to whether to do it will depend on what data the algorithms
are likely to encounter, and that is a statement that holds more generally,
as the next section discusses.

\section{Heuristics}
\label{secheuristics}

A really difficult question in group theory is ``What is a random group?''
Even if we make assumptions about finiteness, and whether we care about
isomorphism types, it is basically impossible to decide in which context one
works. Are we bounding group order or representation degrees?  Are we
picking uniformly from among isomorphism types? Subgroups of $S_n$ or
$GL(n,\mathbb{F})$, and if so for what $n$ and $\mathbb{F}$? Presentations with random relator
strings? But if we cannot answer this question, we also cannot say what a
``typical'' group would be.

This is a question that someone writing {\GAP} code encounters all the
time. Algorithms might be able to take shortcuts if a particular property is
tested for, or it might be possible to trade memory for speed, if not too
much memory would have to be used. The decision on whether to do this will
depend on what objects the algorithm will encounter, and this is a question
every researcher will answer differently.

This uncertainty makes it hard to choose heuristics for how to
proceed within an algorithm. Ultimately, {\GAP} contains many
heuristics that have been driven by the examples the system's authors were
interested in, and different authors would have made different choices. What
seems equally difficult is leaving such a decision to the general user.

Just as it is hard to say when {\GAP} should automatically test for
useful properties, it is unclear
what assumptions we can, or should, make when designing algorithms.

This concerns not only mathematical properties, but also questions of
feasibility. For example, how long should we search for an intermediate
subgroup (if one exists, a long transversal can be
replaced by two shorter ones) when such a subgroup might not exist at all?
When searching for a faithful permutation representation of a factor group,
what degree should we consider as ``good enough'' to stop searching further?

In the end, the heuristics chosen are often influenced by the applications that
motivated the particular implementation, and thus are likely not optimal for
every case. User feedback might make us aware of other relevant cases in
which a particular algorithm does not perform well, and papers describing
the algorithms often come with examples the algorithm's author considered
``standard''. But in general, this is a question that does not seem to have
a good answer.

As far as software is concerned, it will include (sometimes crucial) choices
that were based on the authors' interests, but they can be stretched by
constructive feedback.

\section{Transitioning}
\label{sectrans}

Software is written by humans, and humans have careers (and life spans),
which implies that their involvement with a software system will not be
forever. This poses several questions.

The first of these is: How do we ensure that the production
of useful software is counted, by the community and by hiring and promotion
committees, on an appropriate level, in comparison to the production of
papers? I sadly do not have a good answer here, and while there have been
laudable initiatives, this seems to remain an unsolved problem of
academia.

The second question is how to replace lost capacity and get new contributors involved. The open
source structure of {\GAP} and the ability to contribute code help with
this, but most contributions are short-term, and few contributors
are interested in maintenance tasks. {\GAP} has seen multiple handovers
of responsibility to different groups of people, and some influx of new contributors, though more would be
desirable. I would love to attribute the success of handovers to strategy
and planning. Realistically, it might have been luck, combined with an overall
willingness to put the system's success ahead of personal ambitions.

The third question is what to do with code written in the past by people
who do not maintain it any longer. (This is probably the biggest challenge
to the model of having graduate students contribute to software as part of
their theses.)  We addressed this in part by having others adopt the code,
in part by rewriting. With the benefit of hindsight, I would currently advocate for
rewriting as the ultimately less painful and easier option -- since there is
already working code, this does not necessarily require a redesign from
scratch, and a rewrite can also take care of making the routines more
general or cleaning up a messy interface.
My hope is that in the future, AI systems will be
able to help in improving the documentation and structure of existing code,
making adoption a more plausible and easier choice.

In the end, longevity of a system depends on its community's willingness to
work with the existing code and keep it relevant, at least as much as on
the quality of the code.

\section{Vibe-coding {\GAP}?}
\label{vibecode}

Recent LLM-based AI systems, such as Claude or ChatGPT, have been credited with progress in both programming and
mathematics. The term ``vibe coding''~\cite{vibecoding} has been associated
with this.

It is thus natural to ask what impact such systems will have on the
development of mathematical software. Will users simply use LLMs to write
their mathematical software, rendering existing systems obsolete? The answer
to the last question, I will argue, is ``no'', but we will also see that there
will be a significant impact.

Since my experience has been mostly
with Claude, I will refer to this system, though I would be highly surprised
if other systems behaved fundamentally differently. Experiments in this
section were done in September 2026 with Claude Opus 5.5. The earlier
version mentioned below was from February 2026.

As to the extent of this impact, the answer is simultaneously
``a lot'' and ``not that much'', and the
discrepancy arises from the uneven skill distribution shown by LLMs. The
issue is not so much hallucination as the fact that these systems are good
at some tasks but bad at others. They are like a very eager but sometimes
careless or overconfident graduate student (not that I would ever have
encountered such a student...).

Some of the {\GAP} developers have started to use LLMs (predominantly Claude)
to help with certain tasks, or even to propose code, but any such code
clearly has to undergo a final check by humans before being incorporated
into the system.

I would like to start by looking at tasks that work well. This includes the
category of reviewing and annotating: Adding explanatory comments to code,
documenting existing functions for the manual, renaming variables, or
checking for particular uses of a function all work well and are done much
faster (and, frankly, more completely, given that humans get tired) than a
person could manage. This indeed makes it possible to improve the state of existing
code retroactively, even if the original author has since left academia,
and thus helps with the issues of Section~\ref{sectrans}.

The next level of tasks that work well is producing {\GAP} code from a very
clear code description (e.g. a program in another language). Indeed, my own
Magma-to-{\GAP} translator will likely end up obsolete thanks to the skills of
LLMs.

Claude has broad knowledge of existing functions (its training data
presumably includes
{\GAP}'s source), and it can check the code it writes by running its own instance
of {\GAP} in the cloud.

Such success also extends to routines with a clear structure that mainly rely on a good
understanding of documented information. This could be code to perform a
described set of operations on given input, or code to pre- or post-process
data.

As an example of what is easily possible, I asked Claude to find groups with
isomorphic subgroup lattices:
\begin{quote}
{\it
Write a GAP routine that tests the groups of a given order on whether they have isomorphic subgroup lattices (just the inclusion structure, not subgroup orders)
}
\end{quote}
Claude wrote 131 lines of code (see
\url{https://www.math.colostate.edu/~hulpke/llm/latiso.g}).
It runs through the groups of the given order in the small groups library,
computes their subgroup lattices, translates each lattice into a directed
graph, and leaves the isomorphism test to the {\em Digraphs}
package~\cite{digraphspkg}, which
provides a canonical form that speeds up iterated isomorphism tests.

The routine is reasonably effective, with the runtime dominated by the
subgroup lattice computations. Order 64 takes about 2 minutes (and finds 17
lattices that are shared by 2 or more groups).
For groups whose subgroup lattices are modular,
the results agree with the theory in the
literature~\cite[Chapter 2]{schmidtlattices}.
But that is rather easy code in which all hard work is diverted to
existing {\GAP} functions.
\smallskip

What about a task that I would dread as extremely tedious? {\GAP} has a
function to read in spreadsheets in the \texttt{.csv} format (which is
easily understood from looking at a single example). But the
standard \texttt{.xlsx} format is a compressed archive I did not know
anything about; its documentation on the web~\cite{xlsxformat} is long, and understanding the
format and writing code that parses it would probably take me the
better part of a month.

I asked Claude:

\begin{quote}
{\it
Write a {\GAP} routine, similar to ReadCSV, that can read in an Excel
spreadsheet (\texttt{.xlsx} format) and put the data into a list of records.
}
\end{quote}

Claude ran for about 15 minutes and produced (and
tested) 800 lines of code (see
\url{https://www.math.colostate.edu/~hulpke/llm/readxlsx.g}).
It figured out
that \texttt{unzip} (or \texttt{bsdtar}; the code allows either) unpacks
\texttt{.xlsx} files into
parsable XML. Its code calls this external program and then parses the XML.
Beautiful!

It is important to note, however, that this code requires few design
decisions and little ingenuity, but rather knowledge of the file format, of
the syntax of Unix commands and {\GAP} functions, and willingness to go through
all details to the bitter end.

Notably, there is no consistent programming style. Different sessions will produce
code with different naming and formatting conventions.
\medskip

Perhaps even more impressive is its use as a debugging aid. Claude can look
at code and search for potential errors, in particular if an input that
triggers the error is known. It essentially traces through the code and
identifies where invalid conditions (such as misindexing, or misinterpreting
information, but possibly even some group-theoretic properties) could arise.
It will point out potential locations of errors, but can be wrong. In
particular, it might propose easy fixes that only work around the problem
without touching the underlying issue. Still, this can be an enormous help.

For example, consider the {\GAP} library file \texttt{polyrat.gi}, which implements
polynomial gcd and univariate factorization over the rationals. Most of this
file was written over 25 years ago and has received minimal changes since;
for at least the last decade it could be considered settled code with no
reported errors.

Claude identifies 20 places as problematic. This includes one clear error in
the calculation of bounds for roots of polynomials (caused by misapplication
of absolute values) that will have to be corrected. One recommendation
concerns efficiency. Claude also offers
corrections. Other points fall in the ``cleanup'' versus ``fix''
category.
\smallskip

What much of this has in common is that these are the boring, ``out of
necessity'' parts of programming. They are not the bits that were the reason
why I wanted to program in the first place.
\medskip

So where does an LLM fall short (or
what space remains for humans)?

First, the code it produces simply attempts the task in the most
straightforward way. It does not think about speed, memory use,
feasibility, alternative ways and shortcuts to reduce the search space, or
related questions. (These are the points the previous sections describe
as having been crucial for the system's success.)
It performs the specified task, however good or bad the
task description is.

It does not attempt to change representations, or rephrase the problem, for
efficiency. In short, it treats {\GAP}'s functions as black boxes that answer
questions without ever having to think about
the inside. This is exactly the
``{\em group theoretical fool's paradise where shiny black boxes spit out}
[results]''
that \cite[p.~472]{neubueser93}
warns about.

When it points to potential errors, it might misdiagnose,
or propose unsuitable or inefficient remedies, even if the bug
is in a category that it can often spot. An
example has been an error in calculating the abelianization of a subgroup of
a (complicated) finitely presented group,
\url{https://github.com/gap-system/gap/issues/5609}, which for over a year
was addressed only by a poorly performing temporary workaround, since I was
not able to identify the error.
I asked Claude to help with it, but it only
suggested some ``good-for-you'' cleanup and never identified the true cause,
even though the error (a misindexing in a list, which I eventually found by
going through the code line by line to identify where the example went
wrong) was of a type that Claude had identified successfully in other
places.
\medskip

In short, Claude will write the code you ask it for, but not more. It will not
design sophisticated new algorithms, but
will propose basic calculations that in principle produce the result but might
fail on memory or runtime. If we asked it to expand {\GAP},
unsupervised, it would build many overly generic routines that are limited
in scope.

On the other hand, LLMs will likely allow us in the future to design code
without having to worry about some technical details of language and
libraries. Like previous
transitions, from machine code to high-level languages, or from memory-aware
coding to the use of a memory manager and garbage collector, they will make
programming easier, higher-level, and ultimately more interesting.
In my
youth I spent (too much...) time writing machine code for speed reasons, and
I have never regretted the availability of higher-level languages.
Similarly, future mathematicians will not regret that they can describe
algorithms in more natural language without being tied to a particular
programming language. But users will still need to describe algorithms, not
just tasks, and
the types of objects these algorithms work on. More specifically, they will
have to describe any desired strategies for improved performance.
\medskip

I want to close with an example that might seem to be a perfect outcome,
but leaves some crucial questions open. I took the recent paper
\cite{holttracey26}, which describes an algorithm
for finding a generating set of minimal size for a permutation group. This
paper is an interesting example in that it describes its algorithm
exceptionally well, but never gives pseudocode, and the implementation used
is in Magma (and thus, I believe, not accessible to Claude).

Given only the paper (a PDF, as published), and the prompt to provide a {\GAP}
implementation, Claude produced code (raw version at
\url{https://www.math.colostate.edu/~hulpke/llm/SmallestGeneratingSetNew.g}), and also ran the
smaller tests from the paper for verification.

This code is, in structure and documentation, of high quality. There is,
despite the option to use randomized stabilizer chains, some sluggishness
(which Claude itself acknowledges),
due to repeated generation tests. To make it usable, in particular for the
case of opportune generator reduction that Section~\ref{secgens}
encourages,
more optimization will be needed;
for example, verifying properties of the action (transitivity, primitivity)
before performing element or order tests.

More critically, a first attempt in February 2026
contained a subtle error:
The paper states \cite[p.~132]{holttracey26}
``for each $i=1,\ldots,d$ and $j=1,\ldots,l$ we
test whether $g_1,\ldots,g_{i-1},g_ie_j,g_{i+1},\ldots,g_d$ generate''.
(Here the $g_i$ are chosen representatives of a factor group generating set
and the $e_j$ a generating set of the normal subgroup.)

What is not spelled out (but is implicit in the paper) is that the case
using only the
unmodified generators $g_1,\ldots,g_d$ also has to be tested. (Again it took
significant hand-debugging to find this error, though it is fundamentally of
the misindexing variant.)

To its credit, a newer version of Claude actually spotted this implicit
assumption and took care of it, and this is the code in
\texttt{SmallestGeneratingSetNew.g}.
This clearly shows that careful
inspection of the resulting code, and comparison with the original
description, remain unavoidable.
\bigskip

This brings us back to the question of testing, ensuring correct code (and
making it perform well). Even when using an LLM to write the actual code,
the fundamental tasks of building mathematical software remain the same and
humans will remain crucial.
\medskip

Readers might question how this last statement will hold up in the future.
While predictions are hard, I will go out on a limb and predict
that this observation will persist. The reason for this is that the
(in my view) most interesting algorithms and best implementations in group
theory have always grown with the challenges of specific questions, asked by
humans, which initially often seemed impossible to answer.  The algorithmic
solutions then did not arrive with a snap of the finger but through a,
sometimes arduous, feedback loop of trying code on examples, investigating
where it fails, and trying again, while keeping the interesting examples in
mind. And some of these examples just arise in the process and were never in
the training data.

I also asked this question to Claude. It responded:
\begin{quote}
{\em
There is [...] comparatively little high-quality code for computer algebra
systems [online]. The conventions, idioms and accumulated design decisions
of a long-lived system like GAP or Magma are underrepresented, so the
model’s prior is weakest exactly where the expertise is deepest.
}
\end{quote}
With such an answer, I hope that the existence of this very paper
does not sow the seed for the area's demise.

\section*{AI disclaimer}

In writing this paper I have used Claude as a fancy spell and grammar
checker and occasional thesaurus, and as an incredibly patient proofreader
of earlier versions that demonstrated that English is not my first language.
What it produced were many indications of language problems, which I then
inspected one by one. The writing and the opinions are
entirely mine.

\section*{Acknowledgements}

Many thanks to Alice Niemeyer for extensive feedback on an earlier version
of the paper.
The author's work has been supported in part by
Simons Foundation Grant~852063,
which is gratefully acknowledged.

\end{document}